\documentclass{amsart}

\newtheorem{theorem}{Theorem}[section]
\newtheorem{lemma}[theorem]{Lemma}
\newtheorem{proposition}[theorem]{Proposition}
\newtheorem{corollary}[theorem]{Corollary}

\theoremstyle{definition}

\newtheorem{remarks}[theorem]{Remarks}
\newtheorem*{acknowledgement}{Acknowledgement}

\theoremstyle{remark}

\usepackage{amscd,amssymb}
\usepackage[PostScript=dvips]{diagrams}

\newcommand\mylabel[1]{\label{#1}}

\newcommand{\ZZ}{\mathbb{Z}}

\newcommand{\EE}{\mathbb{E}}
\newcommand{\FF}{\mathbb{F}}
\newcommand{\PP}{\mathbb{P}}
\renewcommand{\AA}{\mathbb{A}}

\newcommand  {\shA}     {\mathcal{A}}

\newcommand  {\shExt}   {\mathcal{E} \!\text{\textit{xt}}}
\newcommand  {\shE}     {\mathcal{E}}

\newcommand  {\shH}     {\mathcal{H}}

\newcommand  {\shI}     {\mathcal{I}}

\newcommand  {\shK}     {\mathcal{K}}

\newcommand  {\shN}     {\mathcal{N}}
\newcommand  {\shL}     {\mathcal{L}}

\newcommand  {\foC}     {\mathfrak{C}}

\newcommand  {\foX}     {\mathfrak{X}}

\newcommand  {\EExt}    {\EE\!\operatorname{xt}}

\newcommand  {\Ext}     {\operatorname{Ext}}

\newcommand  {\Hom}     {\operatorname{Hom}}

\renewcommand  {\k}     {\kappa}

\newcommand  {\lra}     {\longrightarrow}

\newcommand  {\maxid}   {\mathfrak{m}}

\renewcommand{\O}       {\mathcal{O}}

\newcommand  {\Pic}     {\operatorname{Pic}}
\newcommand  {\PGL}     {\operatorname{PGL}}

\newcommand  {\Proj}    {\operatorname{Proj}}

\newcommand  {\quadand} {\quad\text{and}\quad}

\newcommand  {\qM}      {\overline{M}}

\newcommand  {\ra}      {\rightarrow}

\newcommand  {\Sing}    {\operatorname{Sing}}

\newcommand  {\Spec}    {\operatorname{Spec}}
\newcommand  {\Spf}     {\operatorname{Spf}}

\def\mydate{\number\day\space\ifcase\month \or January\or February\or March\or 
April\or May\or June\or July\or
August\or September\or October\or November\or December\fi \space\number\year}

\begin{document}

\title[Triple ramification]
      {Curves with only triple ramification}

\author[Stefan Schroer]{Stefan Schr\"oer}
\address{Mathematische Fakult\"at, Ruhr-Universit\"at, 
         44780 Bochum, Germany}
\curraddr{Mathematisches Institut, Universit\"at zu K\"oln,
Weyertal 86--90, 50931 K\"oln, Germany}
\email{s.schroeer@ruhr-uni-bochum.de}


\dedicatory{Revised version, 2 April 2003}

\begin{abstract}
I show that the set of smooth curves of genus $g\geq 0$ admitting 
a branched covering $X\ra \PP^1$ with only triple
ramification points is of dimension at least
$\max(2g-3,g)$.
In characteristic two,
such curves have   tame rational functions
and an analog of Belyi's Theorem applies to them.
\end{abstract}

\maketitle

\section*{Introduction}

Let $C$ be a smooth proper curve over   an   algebraically
closed ground
field
$k$ of  characteristic $p\geq 0$.
Assuming $p\neq 2$,
Fulton \cite{Fulton 1969} showed that there are
 generically \'etale finite maps $C\ra\PP^1$ such
that all ramification points have index $e=2$.
In this paper, I pose the following question:
\emph{Does there exist a generically \'etale finite map
$C\ra\PP^1$ whose ramification points all have index $e=3$?}

Fried, Klassen, and Kopeliovich \cite{Fried; Klassen; Kopeliovich 2001}
took a first step into this direction. They proved that
all but finitely many complex elliptic curves
admit such a map.
In fact, their proof reveals that for any given genus $g\geq 1$,
the set of Riemann surfaces of genus $g$ admitting such maps is 
at least 1-dimensional.
The arguments, however, are purely topological and
involve homeomorphism spaces, Dehn twists,
and  Teichm\"uller theory.

The main result of this paper is an  improved lower bound
on the dimension via purely algebraic methods. 
We shall prove that the set of points in the moduli space $M_g$ 
whose corresponding curve admits
rational functions with only triple ramification
has dimension  $\geq \max(2g-3,g)$.
Our arguments 
work in all characteristics $p\neq 3$ and rely
on deformation theory and the moduli space of 
stable curves.
The basic idea is to deform
a covering $X_0\ra\PP^1$ where $X_0$ is a curve with  
\emph{cuspidal singularities}, so that each  cuspidal
ramification point breaks up into two regular ramification
points. The stable reduction process involved in this neatly 
explains why we miss $g$ dimensions from the $(3g-3)$-dimensional  
moduli space $M_g$. One might speculate whether   we found
the best lower bound.

Actually, Fulton's result on the existence of maps with ramification
indices
$e=2$ gives some additional information: He proved in \cite{Fulton
1969}, Proposition 8.1, that for a given curve $C$ of genus $g$
in characteristic $p\neq 2$ and $n>g$, there are generically \'etale
maps
$C\ra\PP^1$ of degree $n$ such that all ramification points have index
$e=2$ and that each fiber contains at most one ramification point.
The methods of this paper, however, do not give   much information in
this direction.

My initial motivation to study this problem was \emph{Belyi's Theorem}
\cite{Belyi 1979}.
It states that a compact Riemann surface is defined over
a number field if and only if it admits a finite map
to the Riemann sphere with at most three branch points.
Sa{\"\i}di \cite{Saidi 1997} generalized this to odd characteristics as
follows: An algebraic curve $C$ in characteristic $p\geq 3$ is defined
over a finite field if and only if it admits a \emph{tamely} ramified
morphism
$C\ra\PP^1$ with at most three branch points.
In characteristic $p=2$, the if part holds true, but
the only-if part remains mysterious.
However, a curve $C$ over $\overline{\FF}_2$ admits
a tame function with at most three branch points
if it admits a tame function at all.
In some sense, 
the  result of this paper tells us that the Belyi--Sa{\"\i}di Theorem
is valid in characteristic $p=2$ at least for a $(2g-3)$-dimensional
set.

The question whether a finite morphism $X\ra\PP^1$ whose ramification
points have index $e=3$ exists is also interesting for nonclosed ground
fields. There, however, I showed in \cite{Schroeer 2002} that the
generic curve
$C_\eta$ of
genus $g\geq 3$ in characteristic $p=2$ does not admit such a map. This
relies on   Franchetta's Conjecture, which states that 
$\Pic(C_\eta)=\ZZ K_{C_\eta}$. 
Here the ground field  is the function field $\k(\eta)$ of the moduli
space
$M_g$. Of course, it still might   be true that the desired map exists
over some field extension $\k(\eta)\subset L$.

Here is a plan for the paper.
In Section 1 we study collisions of triple ramification points
in terms of  Weierstrass equations for elliptic curves.
To globalize this, we collect in Section 2 some general results
about deformations of coverings whose fibers are
complete intersections.
We use this to construct effective formal deformation
in Section 3, and explain the resulting increase of
transcendence degree in moduli fields.
Section 4 contains a construction of maps $\PP^1\ra\PP^1$
with only triple coverings so that $\PP^1$ marked
with the ramification points has a large moduli field.
We use this to prove our main result 
in Section 5.
The last section contains some applications regarding Belyi's Theorem
in positive characteristics.

\begin{acknowledgement}
I thank Tam\'as Szamuely, J\'anos Koll\'ar, Fabrizio Catanese, and
Hubert Flenner for helpful discussions,  Le Van Ly for careful
proof-reading, and the referee for his remarks and suggestions,
which helped to clarify the paper. Among other things, he 
suggested to use elliptic curves in Section 
\ref{Collision of triple ramification points} and   moduli spaces  in Section 
\ref{moduli fields}.
\end{acknowledgement}

\section{Collision of triple ramification points}
\label{Collision of triple ramification points}

Fix a 
ground field $k$ of characteristic $p\neq 3$, and let
$h:C\ra D$ be a finite  generically \'etale morphism of proper smooth
curves. For a rational point
$c\in C$ the length $e\geq 1$ of the Artin ring
$\O_{C,c}/\maxid_{D,h(c)}\O_{C,c}$ is called the \emph{ramification
index}. We say that
$h:C\ra D$ has \emph{only triple ramification} if  all ramification
points are rational have ramification index $e=3$.
The key idea of this paper is to collide triple ramification points
in flat families.
We now  explain this by using elliptic curves, where
explicit computations are possible.

Let $A=k[[t]]$ be the formal power series ring in one indeterminate
$t$, and  
$x,y,z$ be homogeneous  coordinates for $\PP^2_A$. Consider the closed
subscheme
$X\subset \PP^2_A$ defined by the Weierstrass equation $x^3=y^2-ty$, or
more precisely by the homogenous equation $x^3=y^2z-tyz^2$.
The generic fiber $X_\eta$ is an elliptic curve with $j$-invariant
$j=0$ over the
field of formal Laurent series $K=k((t))$.
The closed fiber $X_0$ is a rational curve with a cusp located at
$[0,0,1]\in\PP^2_k$.

The diagonal group scheme $G=\mu_{3,A}$ of third roots of unity
acts on
$\PP^2_A$ via the
$\ZZ/3\ZZ$-grading on 
$A[x,y,z]$ given by $\deg(x)=1$ and $\deg(y)=\deg(z)=0$,
as explained in  \cite{SGA 3a}, Expos\'e I, Proposition 4.7.3.
The quotient $\PP^2_A/G$ is the homogeneous spectrum of $A[x^3,y,z]$,
which is isomorphic to a  weighted projective space with weights
$(3,1,1)$. Clearly, the equation $x^3=y^2z-tyz^2$ is homogenous with
respect to the 
$\ZZ/3\ZZ$-grading, so the closed subscheme $X\subset\PP^2_A$ is
$G$-invariant. The corresponding quotient scheme $X/G\subset\PP^2_A/G$
is the homogenous spectrum of
$A[x^3,y,z]/(y^2z-tyz^2-x^3)=A[y,z]$. In particular, we
have
an identification $X/G=\PP^1_A$.

The fixed scheme $X^G $ for the $G$-action on $X$ is the
homogenous spectrum of 
$$
A[x,y,z]/(y^2z-tyz^2-x^3, x)=A[y,z]/(yz(y-tz)).
$$
Hence the generic fiber $X^G_\eta$ comprises the three 
rational points $[0,1,0]$, $[0,0,1]$, and $[0,t,1]$.
In contrast, the closed fiber $X^G_0$ consists of the rational point
$[0,1,0]$ together with  $\Proj k[y,z]/(y^2)$, which is an Artin scheme
of length two around $[0,0,1]$. Intuitively, the two generic fixed
points
$[0,0,1]$ and $[0,t,1]$ collide in the flat family upon specialization.

Let $f:X\ra\PP^1_A$ be the quotient map. Its   generic
fiber
$f_\eta:X_\eta\ra\PP^1_\eta$ is a generically \'etale finite map with
only triple ramification. It has three ramification points,   two of
which collide in  the family. The closed fiber $f_0:X_0\ra\PP^1_k$ is
generically
\'etale as well, but the domain $X_0$ has a cusp resulting from the
collision of triple points.

For later use we take a closer look at the   
complete local rings near 
the cusp $a\in X_0$ and its image $b\in\PP^1_k$.
Clearly, 
$$
\O_{\PP^1_k,b}^\wedge=k[[\frac{y}{z}]]\quadand
\O_{X_0,a}^\wedge=
k[[\frac{x}{z},\frac{y}{z}]]/((\frac{x}{z})^3-(\frac{y}{z})^2).
$$
Let $\tilde{X}_0\ra X_0$ be the normalization map and
$\tilde{a}\in \tilde{X}_0$ be the preimage of $a\in X_0$.
Then there is a uniformizer
$s\in\O_{\tilde{X}_0,\tilde{a}}^\wedge$ such that the inclusion
$\O_{X_0,a}^\wedge\subset\O_{\tilde{X}_0,\tilde{a}}^\wedge$
is nothing but $x/z\mapsto s^2$ and $y/z\mapsto s^3$.
Summing  up:

\begin{proposition}
\mylabel{local cusp}
Inside $\O_{\tilde{X}_0,\tilde{a}}^\wedge=k[[s]]$, the complete local
subalgebra
$\O_{X_0,a}^\wedge$   
 is generated by   $s^2=x/z$ and $s^3=y/z$,
whereas   
$\O_{\PP^1_k,b}^\wedge$ is  generated by $s^3=y/z$.
\end{proposition}

The next task is to find   regular models for the elliptic
curve
$X_\eta$ over $K=k((T))$
whose closed fiber is a reduced divisor with simple normal crossings.
This is indeed possible after replacing $A=k[[t]]$ by
the finite ring extension $A'=k[[t^{1/3}]]$.
Let $X'_\eta$ be the induced
elliptic curve over $K'=k((t^{1/3}))$. The coordinate change 
$x=t^{2/3}\bar{x}$ and $y=t\bar{y}$ shows that 
$\bar{x}^3=\bar{y}^2-\bar{y}$ is a Weierstrass
equation for
$X'_\eta$, so the corresponding constant
elliptic curve over $A'$ is a regular model.

For   later applications, however, we prefer
a regular model $X'\ra\Spec(A')$ such that the projection
$X'_\eta\ra X_\eta$ extends to a morphism $X'\ra X$.
Recall that $X\subset\PP^2_A$ is defined by the homogeneous equation
$x^3=y^2z-tyz^2$. It follows that $\Sing(X)$ consists of a
single point located at $[0,0,1]$, and this singularity
is a rational double
point of type $A_2$. Consider the induced surface $X\otimes A'$.

\begin{proposition}
\mylabel{minimal resolution}
The surface $X\otimes A'$ is normal and   
$\Sing(X\otimes A')$ consists of a single point mapping to the
singular point on $X$. The exceptional divisor $E'\subset X'$ for the
minimal resolution
$X'\ra X\otimes A'$ is an elliptic curve with $j$-invariant $j=0$ and
selfintersection number $-1$. The strict transform $R'\subset X'$
of $X_0$  is a smooth rational curve with selfintersection number $-1$.
The closed fiber $X'_0=E'+R'$ for the projection
$X'\ra\Spec(A') $ is a simple normal
crossing divisor with
$E'\cdot R'=1$.
\end{proposition}

\proof
First note that $X\otimes A'$ satisfies Serre's condition $(S_2)$ by
\cite{EGA IVb}, Proposition 6.4.1. It also satisfies the regularity 
condition $(R_1)$, because   $X\ra\Spec(A)$ is
smooth outside $\Sing(X)$. 
It follows that $X\otimes A'$ is normal and
that
$\Sing(X\otimes A')$ consists of at most one point, which must map to
$\Sing(X)$.

Let $W\ra X$ be the blowing up of the unique singular point
$[0,0,1]\in X$.
A local computation shows that the surface $W$ is regular,
and that the closed fiber $W_0$ is a 
\emph{degeneration of type IV} of the elliptic curve $W_\eta=X_\eta$. 
In other words $W_0=D_1+D_2+D_3$, where  the  $D_i$ are   smooth
rational curves with selfintersection numbers 
$D_i^2=-2$ meeting in a single point. 
Let $V\ra W$ be the blowing up of this intersection point.
Then $V_0=3F+C_1+C_2+C_3$, where $F$ denotes the reduced exceptional divisor
and the $C_i$ are the strict transforms of the $D_i$.
Note that $F^2=-1$ and $C_i^2=-3$.

To get rid of the multiplicity of $F\subset V_0$ we make
a stable reduction process as discussed in \cite{Harris;
Morrison 1998}, Section 3C. Consider the base change $V\otimes A'\ra V$
along the $A$-algebra $A'=k[[t^{1/3}]]$. This is a Kummer covering of
degree three with branch locus  the closed fiber $V_0\subset
V$. We have 
$\O_{V\otimes A'}=\O_V\oplus\O_Vt^{1/3}\oplus\O_Vt^{2/3}$,
where the multiplication law is given by the canonical inclusion
$\O_Vt\subset\O_V$. According to \cite{Esnault; Kahn;
Viehweg 1993}, Proposition 4.3, the normalization $V'\ra V\otimes A'$ is
given by the
$\O_V$-algebra
$$
\O_{V'}=\O_V\oplus\O_V(-C)\oplus\O_V(-2C),
$$
where $C=C_1+C_2+C_3$,  and the 
multiplication law is induced by the composite mapping
$\O_V(-3C)\subset\O_Vt\subset\O_V$.
The projection $h:V'\ra V$ is a Kummer covering of degree three
with branch locus
$C\subset V$. The resulting surface $V'$ is regular because the branch
locus is regular.

Let $F',C_i'\subset V'$ be the reduced preimages of
$F,C_i\subset V$, respectively. We have $h_*(F')=3F$ and $h_*(C'_i)=C_i$
and infer that $V'_0=F'+C_1'+C'_2+C_3'$ is reduced.
By construction, $F'\ra F$ is a Kummer covering of degree three
whose branch locus consists of three reduced points, hence $F'$ 
is an elliptic curve. Since $V'_0$ has arithmetic genus one, the other
components
$C'_i$ are smooth rational curves and the closed fiber
$V'_0$ is a divisor with simple normal crossings. Using
$h^*(F)=F'$ and
$h^*(C_i)=3C_i'$,  we conclude
${F'}^2=-3$ and ${C_i'}^2=-1$ with the projection formula.
After contracting the $(-1)$-curves $C'_i\subset V'$, we obtain 
a relative elliptic curve with generic fiber $V'_\eta=X'_\eta$,
which must be isomorphic to the constant family
with Weierstrass equation $x^3=y^2-y$. In particular, the elliptic
curve $F'$ has $j$-invariant $j=0$.

Let $C_1'\subset V'$ be the reduced strict transform of
$X_0\subset X$, and  $V'\ra X'$ be the contraction of the remaining two
disjoint
$(-1)$-curves
$C'_2\cup C'_3$.  
Then we have an induced map $X'\ra X\otimes A'\ra X$.
The image $E'\subset X'$ of $F'\subset V'$ is an 
elliptic curve with $j$-invariant $j=0$, which is
the exceptional divisor for $X'\ra X\otimes A'$. It follows that
$X\otimes A'$ actually has a singularity and that 
$X'\ra X\otimes A'$ is the minimal resolution of singularities.
The image $R'\subset X'$ of $C_1'$ is a smooth rational curve,
which is also the strict transform of $X_0$.
By construction, $X'_0=E'+R'$ is a simple normal crossing divisor with
intersection numbers  
${E'}^2=-1$, ${R'}^2=-1$, and $E'\cdot R'=1$.
\qed

\begin{remarks}
\mylabel{elliptic singularity}
Laufer showed in \cite{Laufer 1973}, Theorem 4.1 that the formal
isomorphism class of  the singularity on $X\otimes A'$ is uniquely
determined by the $j$-invariant   $j(E)=0$ and the selfintersection
number
$E^2=-1$.
Wagreich observed in \cite{Wagreich 1970}, page 425 that the minimal 
resolution $X'\ra X\otimes A'$ is not realized by  
blowing-up the  reduce singular locus and normalizing.
More generally, Tomari showed in \cite{Tomari 1985}, Theorem 7.4 that
an elliptic Gorenstein surface singularity is resolved by a succession
of blowing ups with reduced centers and normalizations if and
only if the minimal elliptic cycle on the
exceptional divisor has selfintersection $\leq -2$.
\end{remarks}

\section{Deformations for coverings of complete intersection}

Fix a ground field $k$ of arbitrary characteristic $p\geq 0$.
Let $C$ and $D$ be two curves without embedded components,
and $h:C\ra D$ a flat finite morphism that is generically \'etale.
Then $\Omega^1_{C/D}$ is a coherent skyscraper sheaf supported by
the ramification points $x\in C$.
We shall study infinitesimal deformations of $h:C\ra D$.
Let $R$ be a local Artin $k$-algebra with residue field $k$.
A \emph{deformation} of $h$ over $A$ consists of  a   curve $\foC$
flat and of finite type over $A$, a morphism 
$f:\foC\ra D\otimes A$, and  an isomorphism 
$f\otimes_A k\simeq h$.

Suppose $I\subset A$ is an ideal with $I^2=0$, and
$f:\foC\ra D\otimes A/I$ is a deformation of $g$ over $A/I$.
According to \cite{Illusie 1971}, Proposition 2.1.2.3,
the obstruction  for extending it to a deformation
$f':\foC'\ra D\otimes A$ over $A$ lies in the vector space
of hyperextensions
$$
\EExt^2(L^\bullet_{C/D},\O_C\otimes_k I)=
\EExt^2(L^\bullet_{C/D},\O_C)\otimes_k I.
$$
Here $L^\bullet_{C/D}$ is the \emph{cotangent complex} for $h$.
These obstructions vanish under suitable assumptions.
Recall that $h:C\ra D$ is  a 
\emph{morphism locally of complete intersection}
if for all $c\in C$, the Artin
local ring $\O_{C,c}/\maxid_{h(c)}\O_{C,c}$ is the quotient
of some power series algebra $\k(c)[[t_1,\ldots,t_n]]$ by a regular
sequence (\cite{EGA IVd}, Definition 19.3.6).

\begin{proposition}
\mylabel{obstruction}
If the morphism of curves $h:C\ra D$ is locally of complete
intersection, then the group $\EExt^2(L^\bullet_{C/D},\O_C)$ vanishes.
\end{proposition}

\proof
First note that  the cotangent complex
$L^\bullet_{C/D}$ has a very   explicit form in our  situation.
Consider the coherent $\O_D$-module $\shA=h_*(\O_C)$ and the 
projective $D$-scheme
$P=\PP(\shA)$.  Since $h:C\ra D$ is affine, the invertible
$\O_C$-module
$\O_C$ is very ample with respect to $D$, so there is
a closed embedding $C\subset P$.
As $h:C\ra D$ is flat and finite, the $\O_D$-module
$\shA$ is locally free of finite rank, hence the projection $P\ra D$ is
smooth. 
It  follows from \cite{EGA IVd}, Corollary 19.3.5 that the closed
embedding $C\subset P$ is a regular embedding, because
$h:C\ra D$ is a morphism locally of complete intersection.
Hence the conormal sheaf $\shN_{C/P}=\shI/\shI^2$ is   locally free.
By
\cite{Illusie 1971}, Proposition 3.3.6, the cotangent complex
$L_{C/D}^\bullet$ is quasiisomorphic to the  complex concentrated in
degrees $[-1,0]$ given by the canonical map
$\shN_{C/P}\ra\Omega^1_{P/D}|_C$ from the exact sequence
$$
\shN_{C/P}\lra\Omega^1_{P/D}|_C\lra\Omega^1_{C/D}\lra 0.
$$
The map on the left is injective, because it is generically
injective and $\shN_{C/P}$ is torsion free.
It follows that $\shH^{s}(L^\bullet_{C/D})=0$ for $s\neq 0$ and
$\shH^{0}(L^\bullet_{C/D})=\Omega^1_{C/D}$.

Consider the spectral sequence
$\Ext^r(\shH^s(L^\bullet_{C/D}),\O_C)\Rightarrow
\EExt^{r-s}(L^\bullet_{C/D},\O_C)$.
The sheaves $\shH^s(L^\bullet_{C/D})$ are   skyscraper sheaves
supported by the ramification points  $c\in C$. 
We have $\Ext^0(\shH^s(L^\bullet_{C/D}),\O_C)=0$ because $\O_C$ is
torsion free.
It follows that the edge map
$
\Ext^2(\Omega^1_{C/D},\O_C)\ra
\EExt^2(L^\bullet_{C/D},\O_C)
$ 
from the spectral sequence is
surjective.
Next, consider the spectral sequence
$$
H^r(C,\shExt^s(\Omega^1_{C/D},\O_C))\Longrightarrow
\Ext^{r+s}(\Omega^1_{C/D},\O_C).
$$
The sheaf $\shExt^0(\Omega^1_{C/D},\O_C)$ vanishes, because
$\Omega^1_{C/D}$ is torsion and $\O_C$ is torsion free. The group
$H^1(C,\shExt^1(\Omega^1_{C/D},\O_C))$ vanishes because
$\shExt^1(\Omega^1_{C/D},\O_C)$ has 0-dimensional support. The sheaf
$\shExt^2(\Omega^1_{C/D},\O_D)$ vanishes, because the stalks of
$\Omega^1_{C/D}$ have projective dimension $\leq 1$.
We conclude
that
$\Ext^{2}(\Omega^1_{C/D},\O_C)$ and hence
$\EExt^2(L^\bullet_{C/D},\O_C)$ vanishes.
\qed

\medskip 
Suppose again that $I\subset A$ has square zero,
and that $f:\foC\ra D\otimes A/I$ is a deformation over
$A/I$. According to \cite{Illusie 1971}, Proposition 2.1.2.3, the set of
isomorphism classes of deformation   $f':\foC'\ra D\otimes A$
endowed with an isomorphism $f'\otimes A/I\simeq f$ is an
affine space for the vector space of hyperextensions
$$
\EExt^1(L^\bullet_{C/D},\O_C\otimes_k I)=
\EExt^1(L^\bullet_{C/D},\O_C)\otimes_k I.
$$
This group splits up into   local parts:

\begin{proposition}
\mylabel{split}
Let $a_i\in C$ be the ramification
points      and
$b_i=h(a_i)$ the corresponding branch points. Then 
$\EExt^1(L^\bullet_{C/D},\O_C)\simeq
\bigoplus_i
\Ext^1(\Omega^1_{\O_{C,a_i}^\wedge/\O_{D,b_i}^\wedge},
\O^\wedge_{C,a_i})$.
\end{proposition}

\proof
The spectral sequence 
$\Ext^r(\shH^s(L^\bullet_{C/D}),\O_C)\Rightarrow
\EExt^{r-s}(L^\bullet_{C/D},\O_C)$
gives an exact sequence
$$
0\lra \Ext^1(\Omega^1_{C/D},\O_C)\lra
\EExt^1(L^\bullet_{C/D},\O_C)\lra
\Hom(\shH^{-1}(L^\bullet_{C/D}),\O_C).
$$
The term on the right vanishes because
$\shH^{-1}(L^\bullet_{C/D})$ is torsion and $\O_C$ is torsion
free. 
The coherent $\O_C$-module   $\Omega^1_{C/D}$ is a   skyscraper
sheaf supported by the $a_i\in C$, so 
$
\Omega^1_{C/D} =\bigoplus_{i} 
\Omega^1_{\O_{C,a_i}^\wedge/\O_{D,b_i}^\wedge}
$
and the result follows.
\qed

\section{Construction of effective formal deformations}
\mylabel{Construction of effective formal deformations}

We now apply the results of the preceding section
in the following situation. Fix an algebraically closed ground field
$k$ of characteristic $p\neq 3$, and let $\tilde{C}$
be a proper  smooth 
curve of genus $g\geq 0$ over $k$.
Suppose we have a finite generically \'etale   morphism
$\tilde{h}:\tilde{C}\ra \PP^1_k$ with only triple ramification points.

Let $\tilde{a}\in C$ be such a
ramification point. Applying an automorphism of $\PP^1_k$, we may
assume that
$\tilde{a}$ maps to $0\in\PP^1_k$.  Using that $k$ is algebraically
closed and that
$p\neq 3$, one easily sees that there is a uniformizer
$s\in\O_{\tilde{C},\tilde{a}}^\wedge$ such that $s^3$ is a uniformizer
for the subring $\O_{\PP^1,0}^\wedge\subset
\O_{\tilde{C},\tilde{a}}^\wedge$. 

Following Serre's discussion in
\cite{Serre 1975}, Chapter VI, Section 1.3, we now construct a
singular curve $C$ of arithmetic genus $g+1$ with $\tilde{C}$ as
normalization. The underlying topological space of $C$ is the same as
$\tilde{C}$, but we write $a\in C$ for the point corresponding to 
$\tilde{a}\in \tilde{C}$. The
structure sheaf is the sheaf of subalgebras
$\O_{C}\subset\O_{\tilde{C}}$ such that
$\O_{C,c}=\O_{\tilde{C},c}$ for all $c\neq a$.
In contrast, we set 
$\O_{C,a}=k+\maxid_{\tilde{C},\tilde{a}}^2$. 
Intuitively, $C$ is obtained from $\tilde{C}$ by pinching the
first order infinitesimal neighborhood 
of
$\tilde{a}\in\tilde{C}$. The exact sequence
$$
0\lra\O_{C}\lra\O_{\tilde{C}}\oplus\kappa(a)
\lra \O_{\tilde{C},\tilde{a}}/\maxid^2_{\tilde{C},\tilde{a}}\lra 0
$$
gives an exact sequence
$$
0\lra k\lra H^1(C,\O_{C})\lra
H^1(\tilde{C},\O_{\tilde{C}})\lra 0,
$$
and we infer that $C$ has arithmetic genus $h^1(\O_C)=g+1$.
By construction, the canonical bijective morphism  $\tilde{C}\ra C$
is the normalization,  and the image $a\in C$ of
$\tilde{a}\in\tilde{C}$ is a cuspidal singularity with
$\O_{C,a}^\wedge=k[[s^2,s^3]]$.
It follows that  the morphism $\tilde{h}:\tilde{C}\ra\PP^1_k$
induces a finite generically
\'etale  morphism $h:C\ra\PP^1_k$ with
$h(a)=0$.
Obviously, the flat morphism $h:C\ra\PP^1_k$ is locally of complete
intersection.

\begin{proposition}
\mylabel{effective deformation}
Set $A=k[[t]]$. Then there is a flat family 
$X\ra \Spec(A)$ and a finite $A$-morphism $f:X\ra \PP^1_A$
such that the following holds:
\begin{enumerate}
\renewcommand{\labelenumi}{(\roman{enumi})}
\item 
The map on closed fibers $f_0:X_0\ra\PP^1_k$ is isomorphic to
$h:C\ra\PP^1_k$.
\item The formal completion $\O_{X,a}^\wedge$ is isomorphic to 
$k[[x,y,t]]/(y^2-ty-x^3)$
as   algebra over $\O_{\PP^1_A,0}^\wedge=k[[y,t]]$.
\item 
For every closed point $c\in X_0$ with $ c\neq a$, the formal
completion
$\O_{X,c}^\wedge$ is isomorphic to $\O_{X_0,c}^\wedge[[t]]$ as  
algebra over
$\O_{\PP^1_A,f(c)}^\wedge=\O_{\PP^1_k,h(c)}^\wedge[[t]]$.
\item 
The generic fiber $X_\eta$ is a
geometrically connected smooth  curve of genus $g+1$.
\end{enumerate}
\end{proposition}

\proof
First, we shall construct
a formal flat morphism $\foX\ra\Spf(A)$ and a finite formal 
morphism $\foX\ra\PP^1_k\times\Spf(A)$ with 
properties corresponding to (i)--(iii).
Set $A_n=A/(t^{n+1})$. 
Suppose we already have constructed a  flat $A_{n-1}$-scheme $X_{n-1}$
and a morphism $f_{n-1}:X_{n-1}\ra\PP^1_{A_{n-1}}$ with properties as in
(i)--(iii). According to Proposition \ref{obstruction}, there is a flat
$A_n$-scheme
$X'_n$ and a morphism $f'_n:X_n\ra\PP^1_{A_n}$
whose restriction  to $A_{n-1}$ admits an isomorphism
$\varphi'_n:f'_n\otimes A_{n-1}\simeq f_{n-1}$. 
The set of isomorphism classes of such extensions is a torsor
under 
$$
\bigoplus_i
\Ext^1(\Omega^1_{\O_{C,a_i}^\wedge/\O_{\PP^1_k,b_i}^\wedge},
\O^\wedge_{C,a_i})
$$ 
by Proposition \ref{split},
where $a_i\in C$ are the ramification points and $b_i\in \PP^1_k$ 
are the corresponding branch points.
Hence we may choose, for each ramification point $a_i\in C$,
an element 
$\alpha_i\in\Ext^1(\Omega^1_{\O_{C,a_i}^\wedge/\O_{\PP^1_k,b_i}^\wedge},
\O^\wedge_{C,a_i})$ so that $(\alpha_i)$
applied to the isomorphism class of
$(X'_n,f'_n,\varphi'_n)$ gives another  flat $A_n$-scheme $X_n$ together
with a morphism
$f_n:X_n\ra\PP^1_{A_n}$ satisfying our conditions (i)--(iii).
By induction   we construct for all integers $n\geq 0$
morphisms $f_n:X_n\ra\PP^1_{A_n}$ together with
identifications $\varphi_n:f_n\otimes A_{n-1}\simeq f_{n-1}$.
Such a system is nothing but the desired   morphism of
formal schemes.

Being 1-dimensional, the proper scheme $X_0$ admits an ample
invertible sheaf. The obstruction for extending an invertible
sheaf from $X_{n-1}$ to $X_n$ lies in the group $H^2(X_0,\O_{X_0})=0$.
It follows that there is a formal invertible sheaf on $\foX$
that is ample on $X_0$.
By Grothendieck's Algebraization Theorem
(\cite{EGA IIIa}, Theorem 5.4.5), the formal scheme $\foX$
is the formal completion of a projective $A$-scheme $X$.
Moreover (\cite{EGA IIIa}, Theorem 5.4.1),
the   morphism $\foX\ra\PP^1\times\Spf(A)$ of formal schemes comes from
a morphism of schemes $f:X\ra\PP^1_A$.
Properties (i)--(iii) hold because they depend only   on the
underlying   formal schemes.

Concerning the last property (iv),
observe that $C=X_0$ is geometrically integral, because
$\tilde{C}\ra C$ is a birational universal homeomorphism. Then
$X_\eta$ is geometrically integral as well by
\cite{EGA IVc}, Theorem 12.2.1.
Hence $h^0(\O_{X_0})=h^0(\O_{X_\eta})$ and in turn
$g+1=h^1(\O_{X_0})=h^1(\O_{X_\eta})$ by flatness.
\qed

\medskip
Starting with a proper  smooth  $k$-curve
$\tilde{C}$ of genus
$g\geq 1$    endowed with a
finite generically \'etale morphism
$\tilde{C}\ra\PP^1_k$ with only triple ramification, Proposition
\ref{effective deformation} produces a proper
geometrically connected smooth  curve $X_\eta$ over $k((t))$
of genus $g+1$ endowed with a finite generically \'etale morphism
$X_\eta\ra\PP^1_\eta$ with only triple ramification.
The curve $X_\eta$ defines a morphism $\Spec k((T))\ra M_{g+1}$ into the
moduli space of smooth   curves of genus $g+1$. Let $x\in M_{g+1}$ be
the image point. Its residue field $\k(x)$ is called the \emph{moduli
field} for the smooth  curve
$X_\eta$. The crucial observation   is:

\begin{proposition}
\mylabel{not rational}
The moduli field for the smooth   curve $X_\eta$ has
transcendence degree $\geq 1$ over $k$.
\end{proposition}

\proof
Set $A=k[[t]]$ and consider the family $X\ra\Spec(A)$ constructed above.
The generic fiber $X_\eta$ 
defines a rational map $\Spec(A)\dashrightarrow M_{g+1}$.
This rational map is not necessarily everywhere defined, because
$M_{g+1}$ is not proper. However, the moduli space $\qM_{g+1}$ of stable
curves of genus $g+1$ yields a  compactification
$M_{g+1}\subset\qM_{g+1}$, and   the rational map extends
to a morphism $\Spec(A)\ra\qM_{g+1}$.

To prove our assertion it suffices to check that the latter map
does not factor over a closed point. To see this we   make
base change with respect to the ring extension $A'=k[[t^{1/3}]]$.
In light of Proposition \ref{effective deformation} (ii),
the surface $X$ 
at its singular point $a\in X$ is formally isomorphic to the surface
$X$ studied in Section \ref{Collision of triple ramification points}.
It then follows from Proposition \ref{minimal resolution} that
the induced surface $X\otimes A'$ has a unique singularity,
and its minimal resolution $X'\ra X\otimes A'$ yields
a family of stable curves $X'\ra\Spec(A')$.
More precisely, the exceptional curve for the minimal resolution
is an elliptic curve with $j$-invariant $j=0$.
In turn, the image of the  classifying map $\Spec(A')\ra
\qM_{g+1}$ hits the boundary divisor $\qM_{g+1}-M_{g+1}$, whereas the
generic point maps to the interior $M_{g+1}$. It follows that this
morphism does not factor over a closed point, and therefore the moduli
field of $X_\eta$ is not algebraic.
\qed

\section{Moduli fields for pointed rational curves}
\mylabel{moduli fields}

In this section we discuss another method to achieve large 
moduli fields, namely to use pointed rational curves.
Throughout we fix an integer $n\geq 3$. Recall that an 
\emph{$n$-pointed smooth curve of genus zero} over a scheme $S$ is a
smooth proper map $f:X\ra S$ whose fibers are isomorphic to $\PP^1$,
together with $n$ disjoint sections $X_i\subset X$. The following is
well known:

\begin{lemma}
\mylabel{standard form}
In the above situation, there is a unique isomorphism
$X\ra\PP^1_S$ sending $X_1,X_2,X_3\subset X$ to the constant
sections $0,1,\infty\subset\PP^1_S$, respectively.
\end{lemma}

\proof
The $\O_S$-module $\shE=f_*\O_X(X_1)$ is locally free of rank 2, and
$X=\PP(\shE)$. The sections $X_i\subset X$ correspond via
$X_i=\PP(\shL_i)$ to invertible quotients
$\shL_i=\shE/\shK_i$, and the condition of disjointness means
$\shK_i\cap\shK_j=0$ for
$i\neq j$. In turn, the canonical maps $\shE\ra\shL_1\oplus\shL_2$ and
$\shL_1\leftarrow\shK_2\ra\shL_3$ are bijective, so we may assume
$\shL_1=\shL_2$,
$\shE=\shL_1\oplus\shL_1$, and that
$\shK_2\subset\shE$ is the diagonal submodule.
Tensoring   $\shE$ and the $\shL_i$ with $\shL_1^{-1}$, we obtain the
desired isomorphism
$X\ra\PP^1_S$ sending $X_1,X_2,X_3$ to $0,1,\infty$, respectively. This
isomorphism is unique because any automorphism
of $\O_S\oplus\O_S$ fixing the summands and the diagonal submodule is
  multiplication by a scalar.
\qed

\medskip
It follows that the functor sending a scheme $S$ to the
set of isomorphism classes of $n$-pointed smooth curves of genus zero
over $S$ is representable by the scheme
$M_{0,n}=(\PP^1\times\ldots\times\PP^1)-D$. The product has
$n-3$ factors, and $D$ denotes the  closed subset of all  points
$(x_4,\ldots,x_n)$ with $x_i\in\left\{0,1,\infty\right\}$ for some
$4\leq i\leq n$ or  $x_i=x_j$ for some
$4\leq i<j\leq n$.

Fix a ground field $k$.
In light of this explicit nature of $M_{0,n}$, it is easy to compute
moduli fields. Suppose $(\PP^1_K,0,1,\infty,x_4,\ldots,x_n)$ is an
$n$-pointed smooth curve of  genus zero over a field extension
$k\subset K$. Using
the identifications $\AA^1=\PP^1-\left\{\infty\right\}$ and
$\AA^1(K)=K$, we see that the rational points $x_4,\ldots,x_n\in\PP^1_K$
correspond to certain scalars $t_4,\ldots,t_n\in K-\left\{0,1\right\}$.

\begin{proposition}
\mylabel{moduli field}
The moduli field of $(\PP^1_K,0,1,\infty,x_4,\ldots,x_n)$ is
nothing but the subfield $k(t_4,\ldots,t_n)\subset K$.
\end{proposition}

\proof
Let $L\subset K$ be the moduli field in question, and
$(X,x'_1,\ldots,x'_n)$ be a pointed smooth curve of genus zero over $L$
inducing our given pointed curve over $K$.
We may assume $x_1'=0,x_2'=1,x_3'=\infty$ by Lemma \ref{standard form}.
Then the remaining $x_4',\ldots,x_n'$ correspond to scalars
$t'_4,\ldots,t'_n\in L-\left\{0,1\right\}$.
The uniqueness in Lemma \ref{standard form} implies $t_i'=t_i$,
so we have $k(t_4,\ldots,t_n)\subset L$. The reverse inclusion is
obvious.
\qed

\medskip
We now examine the effect of finite morphisms on moduli fields:

\begin{proposition}
\mylabel{transcendence degree}
Let $f:\PP^1_K\ra\PP^1_K$ be a finite morphism 
and $x_1,\ldots,x_n\in\PP^1_K$  be rational points such that the images
$y_1=f(x_i),\ldots,y_n=f(x_n)$ are pairwise different. Then the moduli
fields for the $n$-pointed smooth curves 
$(\PP^1_K;x_1,\ldots,x_n)$ and
$(\PP^1_K;y_1,\ldots,y_n)$ have the same transcendence degree over $k$.
\end{proposition}

\proof
The  finite morphism
$f^{n}:\prod_{i=1}^n\PP^1_K\ra\prod_{i=1}^n\PP^1_K$
induces on the moduli space  a rational map
$\psi:M_{0,n}\dashrightarrow M_{0,n}$ sending an
$n$-pointed smooth curve of genus zero
$(\PP^1,x_1,\ldots,x_n)$ to
$(\PP^1,y_1,\ldots,y_n)$. 
The domain of definition for $\psi$ comprises those
pointed curves for which the image points $y_1,\ldots,y_n$   are
pairwise different. 

The rational map
$\psi:M_{0,n}\dashrightarrow M_{0,n}$ is quasifinite on its domain of
definition, because $f$ is a finite map.
So if $a\in M_{0,n}$ is a point in the domain of definition, and
$b=\psi(a)$ is its image point, then $\k(b)\subset\k(a)$ is a finite
field extension.
This immediately implies our assertion.
\qed

\medskip
We now relate these moduli fields to coverings with
only triple ramification. For this we assume that
our ground field $k$ has characteristic $p\neq 3$.  
Choose an algebraically closed extension field $k\subset K$
of transcendence degree $n-3$, and   algebraically
independent elements $t_4,\ldots,t_n\in K$.
Consider the $n$-pointed smooth curve $(\PP^1_K,y_1,\ldots,y_n)$
with $y_1=0,y_2=1,y_3=\infty$ and such   the
remaining marked points $y_4,\ldots,y_n\in\PP^1_K$ correspond to the
scalars
$t_4,\ldots,t_n\in K$.

\begin{proposition}
\mylabel{covering}
There is a   generically
\'etale finite map
$h:\PP^1_K\ra\PP^1_K$  with only triple ramification such that  the
$y_i\in\PP^1_K$ occur as  branch points. 
If we  choose for each  $y_i\in\PP^1_K$
a ramification point $x_i\in h^{-1}(y_i)$,
then the moduli field for    $(\PP^1_K,x_1,\ldots,x_n)$
has transcendence degree $n-3$.
\end{proposition}

\proof
Consider the polynomial map $r:\PP^1_K\ra\PP^1_K$ given by 
$[z_0,z_1]\mapsto[z_0^3,z_1^3]$. Then
$r$ is a generically
\'etale finite map of degree three with only triple ramification, whose
ramification and branch points are $0,\infty\in\PP^1_K$.
For suitable $\varphi\in\PGL_2(K)$, the composition $\varphi r$
realizes any given pair of rational points $a,b\in\PP^1_K$ as branch
locus.

We now construct the desired map $h$ by induction.
Suppose we already have a generically \'etale finite map
$h_i:\PP^1_K\ra\PP^1_K$ with only triple ramification such that
the $y_j$ occur as branch points for
$1\leq j\leq i$ and $h_i$ is \'etale over $y_j$ for $i+1\leq j\leq n$.
Fix  a rational point $x'\in h_i^{-1}(y_{i+1})$.
Choose $\varphi\in\PGL_2(K)$ so that $x'$ is a branch point
for $\varphi r$, but
that $\varphi r$ is \'etale over 
$h_i^{-1}(y_j)$ for all points $y_j$ with $j\neq i+1$. Then 
$h_{i+1}=h_i\varphi r$ is the desired map.

It remains to check the assertion on   moduli fields.
According to Proposition \ref{moduli field},
the moduli field for $(\PP^1_K,y_1,\ldots,y_n)$ is
$k(t_4,\ldots,t_n)\subset K$, which has transcendence degree $n-3$.
It then follows from Proposition \ref{transcendence degree}
that the moduli field for $(\PP^1_K,x_1,\ldots,x_n)$ has transcendence
degree $n-3$ as well.
\qed

\section{Curves with only triple ramification}

We come to the main result of this paper:

\begin{theorem}
\mylabel{main result}
Let $k$ be a field of characteristic $p\neq 3$.
For each integer $g\geq 0$, there is finitely generated field extension
$k\subset K$ and a smooth  proper curve $C$ of genus $g$ over $K$
with the following properties:
\begin{enumerate}
\renewcommand{\labelenumi}{(\roman{enumi})}
\item The moduli field of 
$C$ has transcendence degree over $k$ at least $\max(2g-3,g)$.
\item There is a finite generically \'etale map
$C\ra\PP^1_K$ with only triple ramification.
\end{enumerate}
\end{theorem}

Let me reformulate this in terms of moduli spaces
over algebraically
closed ground fields.

\begin{corollary}
\mylabel{higher genus}
Suppose $k$ is algebraically closed  of characteristic $p\neq 3$ and
assume $g\geq 2$. Let
$S\subset M_g$ be closure for the  set of all closed points such
that the corresponding curve
$C$ admits a finite generically
\'etale map
$C\ra\PP^1_k$ with only triple ramification.
Then  we have $\dim(S)\geq 2g-3$.
\end{corollary}

\proof
Let $C\ra\PP^1_K$ be as in Theorem \ref{main result},
and choose an integral $k$-scheme $U$ of finite type
whose field of rational functions is $K=\k(U)$.
Shrinking $U$, we may extend $C$ to a smooth relative curve $X\ra U$,
and $C\ra\PP^1_K$ to  a $U$-morphism $f:X\ra\PP^1_U$.
Shrinking further, we may assume that all fibers $f_u:X_u\ra\PP^1_u$
are generically \'etale finite maps with only triple
ramification.

By Chevalley's Theorem, the image $V=\varphi(U)$ of the classifying map
$\varphi:U\ra M_g$ is constructible. Shrinking $U$, we may assume that
$V\subset M_g$ is a subscheme. Since the moduli field of $C$
has transcendence degree $\geq 2g-3$, the dimension of $V$
is at least $2g-3$. For each rational point $\sigma\in V$, the fiber
$\varphi^{-1}(\sigma)\subset U$ contains a rational point because $k$ is
algebraically closed, and the result follows.
\qed

\medskip
We also extend the result on elliptic curves
of Fried, Klassen, and
Kopeliovich \cite{Fried; Klassen; Kopeliovich 2001} 
to all characteristics:

\begin{corollary}
\mylabel{elliptic}
Suppose $k$ is algebraically closed of characteristic $p\neq 3$.
Then for all but finitely many $j$-invariants $j\in k$,
the corresponding elliptic curve
$E$ admits a finite generically \'etale map
$E\ra\PP^1_k$ with only triple ramification.
\end{corollary}

\proof
The argument is as for the preceding corollary, except that one
uses the moduli space $M_{1,1}$ instead of $M_{g}$.
\qed

\medskip\noindent
\emph{Proof of Theorem \ref{main result}.}
First consider the case $g\geq 3$.
We shall construct by induction on $n\geq 0$
a  pointed smooth stable curve
$(C_n,c^n_1,\ldots,c^n_{g-n})$ of genus $n$ with 
$g-n$ marked points over some algebraically closed field
extension $K_n$, so that the moduli field has transcendence
degree
$\geq n+g-3$. Moreover, there will be a generically \'etale finite
map $h_n:C_n\ra\PP^1_{K_n}$ with only triple ramification such that the
$c^n_i$ are ramification points. Induction terminates at $n=g$.

According to Proposition \ref{covering}, the desired
curve exists for $n=0$.
Suppose we already found by induction
$(C_n,c^n_1,\ldots,c^n_{g-n})$ and $h_n:C_n\ra\PP^1_{K_n}$ for some
$n<g$.  The idea now is to trade the last marked point for
a genus increase.
Conforming with the notation in
Section 
\ref{Construction of effective formal deformations},
we set $\tilde{C}=C_n$ and $\tilde{a}=c^n_{g-n}$, and let
$C$ be the corresponding cuspidal curve of genus $n+1$
with normalization $\tilde{C}$.
According to Proposition \ref{effective deformation},
there is an effective deformation $X\ra\Spec(A)$
over $A=K_n[[t]]$ with closed fiber isomorphic to $C$.
Moreover, our given map $C\ra\PP^1_{K_n}$ extends to a family
$f:X\ra\PP^1_A$ whose generic fiber $X_\eta\ra\PP^1_\eta$ is a
generically
\'etale finite map with only triple ramification.
Furthermore, the rational ramification points
$c^n_1,\ldots,c^n_{g-n-1}\in C$ extend to ramification sections over
$A$, and define rational ramification points
$c^{n+1}_1,\ldots,c^{n+1}_{g-(n+1)}\in X_\eta$ in the generic fiber.

Let $X'\ra\Spec(A')$ be the stable reduction over the base change
$A'=k[[t^{1/3}]]$ for $X\ra\Spec(A)$
constructed in the proof for Proposition \ref{not rational}.
Then we have a classifying morphism $\Spec(A')\ra\qM_{n+1,g-(n+1)}$.
The image of the closed point $0\in\Spec(A')$ is a point
$\sigma\in\qM_{n+1,g-(n+1)}$ corresponding to the pointed stable curve  
$(C_n\cup E,c^n_1,\ldots,c^n_{g-(n+1)})$ of genus $n+1$.
Here $E$ is an elliptic curve with $j=0$ as in
the proof for Proposition \ref{not rational}, and 
$C_n\cap E=\left\{c^n_{g-n}\right\}$.
It follows that the residue field $\k(\sigma)$ has transcendence degree
$\geq n+g-3$. As a consequence, the image of the generic point
$\eta\in\Spec(A')$ in $\qM_{n+1,g-(n+1)}$ has   residue field of
transcendence degree
$\geq(n+1)+g-3$.
We now let $K_{n+1}$ be the algebraic closure of $K_n((T))$, set
$C_{n+1}=X_\eta\otimes K_{n+1}$, and choose as marked points
$c^{n+1}_1,\ldots,c^{n+1}_{g-(n+1)}$.
This completes the induction.

In this way we obtain a smooth proper curve $C$ of genus $g$
satisfying properties (i) and (ii). The field extension $k\subset K$ in
the construction is algebraically closed, but it
follows from \cite{EGA IVc}, Theorem 8.8.2  
that the map $C\ra\PP^1_K$ is already defined
over some finitely generated field extension. This finishes
the case $g\geq 3$.

It remains to treat the case $g\leq 2$.
For $g=0$, we simply take $C=\PP^1_k$ and the identity map
$C\ra\PP^1_k$. For $g=1$ or $g=2$, we choose $C_0=\PP^1$,
$c_1^0=0$, $c_2^0=1$, and $c_3^0=\infty$ without
caring for the moduli field, and apply the 
preceding deformation argument
once ore twice, respectively.
\qed

\section{Connections with Belyi's Theorem}

Belyi's Theorem \cite{Belyi 1979} states that
a compact Riemann surface is defined over a number
field if and only if it admits a rational function with
at most three critical values.
Sa{\"\i}di \cite{Saidi 1997} generalized this to odd characteristics
$p\geq 3$. Let me rephrase the part of his result that holds true
for all characteristics:

\begin{proposition}
Let $k$ be an algebraically closed field of characteristic $p>0$.
\begin{enumerate}
\renewcommand{\labelenumi}{(\roman{enumi})}
\item A smooth proper curve $C$ over $k$ is defined over a finite
field if there is a
finite map $C\ra\PP^1_k$ with only tame ramification and
at most three branch points.
\item A smooth proper curve $C$ over $\overline{\FF}_p$ admits
a 
finite map $h:C\ra\PP^1_{\overline{\FF}_p}$ with only tame ramification
and at most three branch points if
there is at least one finite  map $g:C\ra\PP^1_{\overline{\FF}_p}$
with only tame ramification.
\end{enumerate}
\end{proposition}

\proof
For convenience, I recall Sa{\"\i}di's argument:
The first statement follows from
Grothendieck's theory of the tame fundamental group
(compare \cite{Orgogozo; Vidal 2000}, Theorem 6.1).
For the second statement, let
$h_n:\PP^1_{\overline{\FF}_p}\ra\PP^1_{\overline{\FF}_p}$ be the
polynomial map
$[z_0,z_1]\mapsto [z_0^{p^n-1},z_1^{p^n-1}]$. Then $h=h_n\circ g$ is
the desired map for some $n$ sufficiently large, as explained in
\cite{Saidi 1997}, Theorem 5.6.
\qed

\medskip
In characteristic $p\geq 3$, tame functions
$g:C\ra\PP^1_{\overline{\FF}_p}$ as in (ii)  exist by
\cite{Fulton 1969}, Proposition 8.1. In characteristic $p=2$,  
Corollary \ref{higher genus} tell us that this holds
true for a set of of curves of dimension at least $2g-3$.


\end{document}